\documentstyle[12pt]{amsart}
\numberwithin{equation}{section}
\begin{document}

%



\setlength{\itemsep}{0in}\newcommand{\lab}{\label}
\newcommand{\labeq}[1]{  \be \label{#1}  }
\newcommand{\labea}[1]{  \bea \label{#1}  }
\newcommand{\ben}{\begin{enumerate}}
\newcommand{\een}{\end{enumerate}}
\newcommand{\bm}{\boldmath}
\newcommand{\Bm}{\Boldmath}
\newcommand{\bea}{\begin{eqnarray}}
\newcommand{\ba}{\begin{array}}
\newcommand{\bean}{\begin{eqnarray*}}
\newcommand{\ea}{\end{array}}
\newcommand{\eea}{\end{eqnarray}}
\newcommand{\eean}{\end{eqnarray*}}
\newcommand{\bthm}{\begin{thm}}
\newcommand{\ethm}{\end{thm}}
\newcommand{\blem}{\begin{lem}}
\newcommand{\elem}{\end{lem}}
\newcommand{\bprop}{\begin{prop}}
\newcommand{\eprop}{\end{prop}}
\newcommand{\bcor}{\begin{cor}}
\newcommand{\ecor}{\end{cor}}
\newcommand{\bdfn}{\begin{dfn}}
\newcommand{\edfn}{\end{dfn}}
\newcommand{\brem}{\begin{rem}}
\newcommand{\erem}{\end{rem}}
\newcommand{\lb}{\linebreak}
\newcommand{\nlb}{\nolinebreak}
\newcommand{\nl}{\newline}
\newcommand{\hs}{\hspace}
\newcommand{\vs}{\vspace}
\alph{enumii} \roman{enumiii}
\newtheorem{thm}{Theorem}[section]
\newtheorem{prop}[thm]{Proposition}
\newtheorem{lem}[thm]{Lemma}
\newtheorem{sublem}[thm]{Sublemma}
\newtheorem{cor}[thm]{Corollary}
\newtheorem{dfn}[thm]{Definition}
\newtheorem{rem}[thm]{Remark}

\def\endpf{\qed}
\def\ms{\medskip}
\def\N{I\!\!N}                \def\Z{Z\!\!\!\!Z}      \def\R{I\!\!R}
\def\C{C\!\!\!\!I}            \def\T{T\!\!\!\!T}      \def\oc{\overline \C}
\def\Q{Q\!\!\!\!|}
\def\1{1\!\!1}

\def\and{\text{ and }}        \def\for{\text{ for }}
\def\tif{\text{ if }}         \def\then{\text{ then }}

\def\Cap{\text{Cap}}          \def\Con{\text{Con}}
\def\Comp{\text{Comp}}        \def\diam{\text{\rm {diam}}}
\def\dist{\text{{\rm dist}}}  \def\Dist{\text{{\rm Dist}}}
\def\Crit{\text{Crit}}
\def\Sing{\text{Sing}}        \def\conv{\text{{\rm conv}}}

\def\Fin{{\cal F}in}
\def\F{{\Cal F}}
\def\h{{\text h}}
\def\hmu{\h_\mu}           \def\htop{{\text h_{\text{top}}}}

\def\H{\text{{\rm H}}}     \def\HD{\text{{\rm HD}}}   \def\DD{\text{DD}}
\def\BD{\text{{\rm BD}}}         \def\PD{\text{PD}}
\def\re{\text{{\rm Re}}}    \def\im{\text{{\rm Im}}}
\def\Int{\text{{\rm Int}}} \def\ep{\text{e}}
\def\CD{\text{CD}}         \def\P{\text{{\rm P}}}     \def\Id{\text{{\rm Id}}}
\def\Hyp{\text{{\rm Hyp}}}
\def\A{\Cal A}             \def\Ba{\Cal B}       \def\Ca{\Cal C}
\def\Ha{\Cal H}
\def\L{{\cal L}}             \def\M{\Cal M}        \def\Pa{\Cal P}
\def\U{\Cal U}             \def\V{\Cal V}
\def\W{\Cal W}             \def\Log{\text{{\rm Log}}}

\def\a{\alpha}                \def\b{\beta}             \def\d{\delta}
\def\De{\Delta}               \def\e{\varepsilon}          \def\f{\phi}
\def\g{\gamma}                \def\Ga{\Gamma}           \def\l{\lambda}
\def\La{\Lambda}              \def\om{\omega}           \def\Om{\Omega}
\def\Sg{\Sigma}               \def\sg{\sigma}
\def\Th{\Theta}               \def\th{\theta}           \def\vth{\vartheta}
\def\ka{\kappa}
\def\Ka{\Kappa}

\def\bi{\bigcap}              \def\bu{\bigcup}
\def\({\bigl(}                \def\){\bigr)}
\def\lt{\left}                \def\rt{\right}
\def\bv{\bigvee}
\def\ld{\ldots\!}               \def\bd{\partial}         \def\^{\tilde}
\def\club{\hfill{$\clubsuit$}}\def\proot{\root p\of}

\def\tm{\widetilde{\mu}}
\def\tn{\widetilde{\nu}}
\def\es{\emptyset}            \def\sms{\setminus}
\def\sbt{\subset}             \def\spt{\supset}

\def\gek{\succeq}             \def\lek{\preceq}
\def\eqv{\Leftrightarrow}     \def\llr{\Longleftrightarrow}
\def\lr{\Longrightarrow}      \def\imp{\Rightarrow}
\def\comp{\asymp}
\def\upto{\nearrow}           \def\downto{\searrow}
\def\sp{\medskip}             \def\fr{\noindent}        \def\nl{\newline}

\def\ov{\overline}            \def\un{\underline}

\def\ess{{\rm ess}}           \def\sd{\bigtriangleup}
\def\ni{\noindent}
\def\cl{\text{cl}}
\def\om{\omega}
\def\bt{{\bf t}}
\def\Bu{{\bf u}}
\def\tr{t}
\def\bo{{\bf 0}}
\def\nut{\nu_\lla^{\scriptscriptstyle 1/2}}
\def\arg{\text{arg}}
\def\Arg{\text{Arg}}
\def\re{\text{{\rm Re}}}
\def\gr{\nabla}
\def\endpf{{${\mathbin{\hbox{\vrule height 6.8pt depth0pt width6.8pt  }}}$}}
\def\Fa{\cal F}
\def\Gal{\cal G}
\def\1{1\!\!1}
\def\D{{I\!\!D}}

\title{The Equilibrium States \\ for \\
Semigroups of Rational Maps}

\thanks{\noindent Date: August 19, 2008. Published in Monatsh. Math., 156 (2009), no. 4, 
371--390.  
The research of the second author was supported in part by  
the NSF Grant DMS 0400481. The first author thanks University of 
North Texas for kind hospitality, during his stay there.}

\

\author{Hiroki Sumi and Mariusz Urba\'nski}
\address{\ni Hiroki Sumi; 
Department of Mathematics,
Graduate School of Science,
Osaka University, 
1-1 Machikaneyama,
Toyonaka,
Osaka, 560-0043, 
Japan
\newline {\tt sumi@@math.sci.osaka-u.ac.jp,
http://www.math.sci.osaka-u.ac.jp/$\sim$sumi/}}

\

\address{\ni Mariusz Urba\'nski; Department of
Mathematics, University of North Texas, P.O. Box 311430, Denton
TX 76203-1430, USA
\newline {\tt urbanski@@unt.edu,
http://www.math.unt.edu/$\sim$urbanski}}

\begin{abstract}
We consider the dynamics of skew product maps associated with 
finitely generated semigroups of rational maps on the Riemann sphere. 
We show that under some conditions on the dynamics and the potential 
function $\psi $, there exists a unique equilibrium 
state for $\psi $ and a unique $\exp(\P(\psi )-\psi )$-conformal 
measure, where $\P(\psi )$ denotes the topological pressure of $\psi$. 
\end{abstract}
\maketitle
\section{{\bf Introduction}}

In this paper, we frequently use the notation from \cite{hiroki1}. 
A ``rational semigroup" $G$ is a semigroup generated by a family of non-constant 
rational maps $g:\oc \rightarrow \oc$,\ where $\oc $ denotes the 
Riemann sphere,\ with the semigroup operation being functional 
composition. For a rational semigroup $G$,\ we set 
$F(G):=\{ z\in \oc \mid G \mbox{ is normal in a neighborhood of } z\} 
$ and $J(G):=\oc \setminus F(G).$ 
$F(G)$ is called the Fatou set of $G$ and $J(G)$ is called the 
Julia set of $G.$ 
 If $G$ is generated by a family $\{ f_{i}\} _{i}$,\ 
 then we write $G=\langle f_{1},f_{2},\ldots\! \rangle .$ 

The research on the dynamics of rational semigroups was initiated by 
Hinkkanen and Martin (\cite{HM}), who were interested 
in the role of the dynamics of polynomial 
semigroups while studying various one-complex-dimensional moduli spaces 
for discrete groups, and by 
F. Ren's group (\cite{ZR}), who studied such semigroups 
from the perspective of random complex dynamics.  
For further studies on the dynamics of rational semigroups, 
see \cite{hiroki1, hiroki2, hiroki3, hiroki4, 
hirokidc, hiroki5, StSu, SU}.

The theory of the dynamics of rational semigroups is  
 deeply related to that of the fractal geometry. In fact, 
 if $G=\langle f_{1},\ldots\!\!,f_{s}\rangle $ is a finitely generated 
 rational semigroup, then the Julia set $J(G)$ of $G$ has the 
 ``backward self-similarity'', i.e., 
 $$J(G)=f_{1}^{-1}(J(G))\cup \cdots \cup f_{s}^{-1}(J(G))$$ 
 (See \cite{hiroki1}). Hence, 
 the behavior of the (backward) dynamics of finitely generated 
 rational semigroups can be regarded 
 as that of the ``backward iterated function systems.''
For example, the Sierpi\'{n}ski gasket can be regarded as 
the Julia set of a rational semigroup.

  The research 
on the  
dynamics of rational semigroups is also directly related to 
that 
on the 
random dynamics of holomorphic maps. 
The first 
study 
on the 
random dynamics of holomorphic maps 
was 
by Fornaess and Sibony (\cite{FS}), and 
much research has followed. 
(See \cite{Br, Bu1, Bu2, 
BBR}.) In fact, it is very natural 
and important to combine both the theory 
of rational semigroups and that of random 
complex dynamics (see \cite{hiroki1, hiroki2, 
hiroki4, hirokidc, hiroki5}).

 We use throughout spherical derivatives
and, for each meromorphic function $\varphi $, 
we denote by $|\varphi '(z)|$ the norm of the derivative 
with respect to the spherical metric.
We denote by 
$CV(\varphi )$ the set of critical values of $\varphi .$ 
The symbol $A$ is used to denote the spherical area. 
We set $\ka=A(B(0,1))/4$. Obviously, 
there exists a constant $C_{sa}\geq 1$ such that 
for each $0<R\leq \mbox{diam} ({\oc })/2$,  
$C_{sa}^{-1}\leq A(B(z,R))/\ka(2R)^2\leq C_{sa}.$

  Let $G=\langle f_{1},\ldots\! ,f_{s}\rangle $ be a 
finitely generated rational semigroup. 
Then,\ we use the following notation. 
Let $\Sigma _{s}:=\{ 1,\ldots\! ,s\} ^{\Bbb{N}}$ be the 
space of one-sided sequences of $s$-symbols endowed with the 
product topology. This is a compact metric space. 
Let 
$f:\Sg_{s}\times \oc \rightarrow \Sg_{s}\times \oc $ 
be the skew product map associated with $\{ f_{1},\ldots\! ,f_{s}\} $
given by the formula,\ $f(\om,z)=(\sg (\om ),\ f_{\om_{1}}(z))$,\ where 
$(\om,z)\in \Sg _{s}\times \oc,\ \om=(\om_{1},\om_{2},\ldots\! ),$ and 
$\sg :\Sigma _{s}\rightarrow \Sg _{s}$ denotes the shift map.
We denote by $\pi _{1}:\Sigma _{s}\times \oc \rightarrow \Sigma _{s}$ 
the projection onto $\Sigma _{s}$ and 
$\pi _{2}:\Sigma _{s}\times \oc \rightarrow \oc $ the 
projection onto $\oc $. That is, 
$\pi _{1}(\om ,z)=\om $ and $\pi _{2}(\om ,z)=z.$  
Under the canonical identification $\pi _{1}^{-1}\{ \om \} \cong \oc $, 
each fiber $\pi _{1}^{-1}\{ \om \} $ is a Riemann surface which 
is isomorphic to $\oc .$ 

 Let Crit$(f):=
\bigcup_{\om \in \Sigma _{s}}\{ v\in \pi _{1}^{-1}\{ \om \} \mid 
v \mbox{ is a critical point of }f
|_{\pi _{1}^{-1}\{ \om \} }\rightarrow \pi _{1}^{-1}\{ \sigma (\om )\} \} 
\ (\subset \Sigma _{s}\times \oc )  
$ be the set of critical points of $f.$  
For each $n\in \Bbb{N} $ and $(\om ,z)\in 
\Sigma _{s}\times \oc $, we set $(f^{n})'(\om ,z):= (f_{\om_{n}}\circ
\cdots \circ f_{\om _{1}})'(z).$  

For each $\om\in \Sg_s$ we define 
$$
J_{\om }:=\{ z\in \oc \mid 
 \{ f_{\om _{n}}\circ \cdots \circ f_{\om _{1}}\} _{n\in \Bbb{N}} \mbox{ is 
 not normal in any neighborhood of } z\} 
$$
 and we then set 
$$
J(f):= \overline{\cup _{w\in \Sigma _{s}}\{ \om \} \times J_{\om } },
$$
 where the closure is taken in the product space 
 $\Sigma _{s}\times \oc .$ By definition,\ 
 $J(f)$ is compact. Furthermore,\ by Proposition 3.2 in \cite{hiroki1},\ 
 $J(f)$ is completely invariant under $f$,\ 
 $f$ is an open map on $J(f)$,\ 
 $(f,J(f))$ is topologically exact under a mild condition,\ and $J(f)$ is equal to the closure of 
 the set of repelling periodic points of 
 $f$ provided that $\sharp J(G)\geq 3$,\ where we say that a periodic point $(\om ,z)$ of $f$ with 
 period $n$ is repelling if $|(f^{n})'(\om ,z)|>1.$      
 Furthermore,\ $\pi _{2}(J(f))=J(G).$ 
We set 
$$
e_{j}:=\deg (f_{j}) \  \text{ and }  \  d:=\max_{j=1,\ldots\!
  ,s}(2e_{j}-2).
$$
Throughout the paper, we assume the following.
\begin{itemize}
\item[(E1)] There exists an element 
$g\in G$ with $\deg (g)\geq 2.$ Furthermore,\ for the semigroup $H:=\{
h^{-1}\mid h\in \mbox{Aut} \oc \cap G\}$,\ we have 
$J(G)\subset F(H).$ (If $H$ is empty, we put $F(H):=\oc .$)
\item[(E2)]$\sharp $(Crit $(f)\cap J(f))<\infty $
\item[(E3)]There is no super attracting cycle of $f$ in $J(f).$
That is, if $(\om ,z)\in J(f)$ is a periodic point with period $n$,\ 
then $(f^{n})'(\om ,z)\neq 0.$  
\end{itemize}

\ 

\ni Any finitely generated rational semigroup 
$G=\langle f_{1},\ldots\! ,f_{s}\rangle $ satisfying all the conditions 
(E1)-(E3) will be called an E-semigroup of rational maps.\\

\noindent {\bf Examples.} 
\begin{itemize}
\item 
If $g$ is a rational map with $\deg (g)\geq 2$,\ 
then $\langle g\rangle $ is an E-semigroup of rational maps.
\item Let $G=\langle f_{1},\ldots\! ,f_{s}\rangle $ be a 
finitely generated rational semigroup such that 
 $\deg (f_{j})\geq 2$ for each $j= 1,\ldots ,s.$  
If CV$(f_{j})\cap J(G)=\emptyset $ for each $j=1,\ldots\! ,s$,\ 
then $G$ is an E-semigroup of rational maps. 
\item Let $g_{1}$ and $g_{2}$ be two polynomials with 
$\deg (g_{j})\geq 2\ (j=1,2)$. Suppose that 
$$
\overline{\cup _{n=0}^{\infty }g_{1}^{n}(CV(g_{1}))} 
\subset A_{\infty }(g_{2}) \  \text{  and } \
\overline{\cup _{n=0}^{\infty }g_{2}^{n}(CV(g_{2}))}\subset A_{\infty
}(g_{1}),
$$
where $A_{\infty }(\cdot )$ denotes the basin of infinity. Let 
$m\in \Bbb{N}$ be a sufficiently large number so that 
$$
g_{1}^{m}(\overline{\cup _{n=0}^{\infty }g_{2}^{n}(CV(g_{2}))})\subset
U_{\infty } \
\text{ and } \ 
g_{2}^{m}(\overline{\cup _{n=0}^{\infty }g_{1}^{n}(CV(g_{1}))})\subset U_{\infty },  
$$
where $U_{\infty }$ denotes the connected component of 
$F(\langle g_{1},g_{2}\rangle )$ with $\infty \in U_{\infty }.$ 
Then,\ the semigroup $\langle g_{1}^{m},g_{2}^{m}\rangle $ is an E-semigroup 
of rational maps.  
\end{itemize}
The dynamics of the skew product map $f:\Sigma _{s}\times \oc 
\rightarrow \Sigma _{s}\times \oc $ is directly related to 
the dynamics of the semigroup $G=\langle f_{1},\ldots\! 
,f_{s}\rangle .$ 
For example, we use the dynamics of $f$ to analyze the 
dimension of Julia set $J(G)$ of $G$ (see 
\cite{hiroki1, hiroki2, hiroki3,  
hiroki4, SU}). 

Our main concern in this paper is the existence and uniqueness of
equilibrium states for all $E$-semigroups of rational maps and a
large class of H\"older continuous potentials. We do not assume any 
kind of expanding property. If $T:X\to X$ is a continuous mapping of a
compact metric space $X$, then one can define the topological pressure
$\P(g)$ for every real-valued continuous function $g$ on $X$ (see
\cite{wa} for instance) in purely topological terms. The link with
measurable dynamics is given by the Variational Principle (see
\cite{wa} again) saying that the topological pressure $\P(g)$ is equal
to the supremum of all numbers $\hmu(T)+\int gd\mu$, $\mu$ being Borel
probability $T$-invariant measures on $X$. Any invariant measure $\mu$
satisfying equality
$$
\P(g)=\hmu(T)+\int gd\mu
$$
is called an equilibrium state. All 
equilibrium states 
have a profound
physical meaning (see \cite{ruelle}), their existence as well as
uniqueness is of primary importance and is in general not clear at
all. A general sufficient condition for the existence is expansiveness,
a weaker one, holding for all rational maps on the sphere (see 
\cite{lyubich}), is asymptotic $h$-expansiveness. This covers the case of
hyperbolic maps studied in \cite{hiroki3}.
Our goal in this paper is to do both, existence and uniqueness of
equilibrium states, for E-semigroup of rational maps and a
large class of H\"older continuous potentials. Thus, the main
purpose of this paper is to prove the following result.

\

\bthm\lab{mainth}({\bf Main result}) 
Let $G=\langle f_{1},\ldots\! ,f_{s}\rangle $ be an E-semigroup 
of rational maps. 
Let $f:\Sigma _{s}\times \oc \rightarrow \Sigma _{s}\times \oc $ be 
the skew product map associated with 
$\{ f_{1},\ldots\! ,f_{s}\} .$ 
Let $\psi :J(f)\rightarrow \Bbb{R} $ be a H\"{o}lder continuous 
function. Moreover, let $P^{p}(\psi )$ be the 
pointwise pressure of $\psi $ 
with respect to the dynamics of $f:J(f)\rightarrow J(f)$ 
(see the definition (\ref{PxPp})). Suppose that 
$\P^{p}(\psi )>\sup (\psi )+\log s.$ 
Then, all of the following statements hold. 
\begin{enumerate}
\item 
Regarding the dynamics of 
$f:J(f)\rightarrow J(f)$, there exists a unique 
equilibrium state for $\psi $ and a unique 
$\exp (\P(\psi )-\psi )$-conformal measure in the sense of 
\cite{duecm}, where $\P(\psi )$ denotes the pressure of 
$(f|_{J(f)}, \psi ).$  
\item Moreover, 
$\P(\psi )=\P^{p}(\psi )$ and 
for each point $x\in J(f)$, 
we have that $\frac{1}{n}\log \L_{\psi }^{n}\1 (x)
\rightarrow \P(\psi )=\P^{p}(\psi )$ as $n\rightarrow \infty $, 
where $\L _{\psi }$ denotes the Perron-Frobenius 
operator associated with the potential $\psi $
 (see {\em (\ref{pfopdf})}) and $\1 \equiv 1.$ 
\end{enumerate}
\ethm

\

\begin{rem}
If $\psi :J(f)\rightarrow \Bbb{R}$ is a H\"{o}lder continuous 
function such that $\sup (\psi )-\inf (\psi )<\log (\sum _{j=1}^{s}e_{j})
-\log s$, then $\P^{p}(\psi )>\sup (\psi )+\log s.$ 
\end{rem}

\

\begin{rem}
Under a very mild condition, 
the topological entropy $h(f)$ of $f: \Sigma_{s}\times \oc \rightarrow 
\Sigma_{s}
\times \oc $ is equal to $\log (\sum _{j=1}^{s}e_{j})$  and 
there exists a unique maximal entropy measure for 
$f: \Sigma _{s}\times \oc \rightarrow \Sigma _{s}\times \oc $ (See \cite{hiroki1}). 
\end{rem}

\

\ni The proof of the main result is given in the following several sections.
In this proof we utilize some arguments from \cite{duecm} to show the 
existence of a conformal measure. Then, using the normality of 
a family of inverse branches of elements of the semigroup 
(see section~\ref{Distth}), we analyze the Perron-Frobenius operator in detail. 
Developing the techniques worked out in
\cite{dues} and \cite{hiroki1}, we show the existence of an equilibrium 
state and the uniqueness of a conformal measure as well as 
an equilibrium state.  
 
\section{Distortion Theorems}
\label{Distth}
\ni Let us recall the following well-known version of Koebe's Distortion
Theorem concerning spherical derivatives.

\

\bthm\lab{t1042805}
For every $u\in(0,\diam(\oc)/2)$ there exists a function 
$k_u:[0,1)\to (0,+\infty)$, continuous at $0$, with $k_u(0)=1$ and the
following property. If $\xi\in\oc$, $R>0$, and $H:B(\xi,R)\to\oc$ is a
meromorphic univalent function such that $\oc\sms H(B(\xi,R))$ contains 
a ball of radius $u$, then for every $t\in[0,1)$ and all 
$z,w\in B(\xi,tR)$,
$$
k_u^{-1}(t)\le{|H'(w)|\over |H'(z)|}\le k_u(t).
$$
\ethm

\

\ni As an immediate consequence of this theorem,\ combined with 
Lemma 4.5 in \cite{hiroki1},\ we get the following.

\
\blem\lab{l3042805}
Let $G=\langle f_{1},\ldots\! ,f_{m}\rangle $ be a finitely generated 
rational semigroup satisfying the condition (E1). Then,\ 
there exists a number $R_{0}>0$ and a 
function $k_{G}:[0,1)\times (0,R_{0}]\rightarrow [1,\infty )$ 
such that 
for each $x\in J(G),0<R\leq R_{0}$ and $0\leq t<1$,\ the 
family 
${\cal F}_{x,R}:= 
\{ \varphi :B(x,R)\rightarrow \oc \mid \varphi \mbox{ is a well-defined 
inverse branch of } h,\ h\in G\} $ satisfies that 
for each $H\in {\cal F}_{x,R}$,\ $w,z\in B(x,tR)$,\ 
we have 
$$ k_{G}(t,R)^{-1}\leq \frac{|H'(w)|}{|H'(z)|}\leq k_{G}(t,R).$$  

\elem
{\sl Proof.}
By Lemma 4.5 in \cite{hiroki1}, 
we have that there exists a number $R_{0}>0$ such that 
for each $x\in J(G)$ and each $0<R\leq R_{0}$, the family 
${\cal F}_{x,R}$ is normal in $B(x,R).$ 
Now, suppose the statement of our lemma is false. 
Then, there exist a $0<t<1$, a sequence 
$(x_{n})$ in $J(G)$, a sequence $(w_{n})$ in $\oc $, 
a sequence $(z_{n})$ in $\oc $, 
a number $0<R\leq R_{0}$, 
and a sequence $(\varphi _{n})$ of 
meromorphic functions, such that for each $n\in \Bbb{N}$, 
we have 
$w_{n},z_{n}\in B(x_{n},tR),\ 
\varphi _{n}\in {\cal F}_{x_{n},R}$, and 
$\frac{|\varphi _{n}'(w_{n})|}{|\varphi _{n}'(z_{n})|}\geq n.$ 
We may assume $x_{n}\rightarrow x_{\infty }\in J(G),\ 
w_{n}\rightarrow w_{\infty }\in \oc,\ $ and 
$z_{n}\rightarrow z_{\infty }\in \oc .$ 
%
Let $u$ be a number with 
 $t<u<1.$ Then there exists an $n_{0}\in \Bbb{N}$ such that 
 each $\varphi _{n}\ (n\geq n_{0})$ is defined on $B(x_{\infty },
 uR)$ and the family 
 $\{ \varphi _{n}\} _{n\geq n_{0}}$ is normal 
 in $B(x_{\infty },uR).$ 
 Then we may assume that $\varphi _{n}$ tends to a meromorphic 
 function $\varphi _{\infty }:B(x_{\infty },uR)\rightarrow 
 \oc $ as $n\rightarrow \infty $,\ uniformly on  
 $B(x_{\infty }, uR).$   
Since each $\varphi _{n}$ is injective, we have that 
$\varphi _{\infty }:B(x_{\infty },uR)\rightarrow \oc $ is 
either injective or constant. Hence, 
$\varphi _{\infty }(B(x_{\infty },uR))\neq \oc .$ 
Let $v$ be a number with $t<v<u.$ 
Then, there exist a number $0<s<1$ and a point $a\in \oc $ such that 
for each large $n\geq n_{0}$, 
$$B(a,s)\subset \oc \setminus \varphi _{n}(B(x_{\infty },vR)).$$ 
Then, by Theorem~\ref{t1042805}, it causes a contradiction.
We are done.
\endpf

\ 

\noindent {\bf Notation.} Throughout the rest of the paper,\ 
we set $K_{G}(R):=k_{G}(\frac{1}{2},R).$ 


\section{Inverse Branches}

We set 
$\Sigma _{s}^{\ast }:= \cup _{j=1}^{\infty }\{ 1,\ldots\! ,s\} ^{j}$ 
(disjoint union). For each $\om \in \Sigma _{s}^{\ast }$, 
we set $| \om | =j$ if $\om \in \{ 1,\ldots\! ,s\} ^{j}.$ 
By $\hat\Sg_s^+$ we
denote the space dual to 
$\Sg_{s}$, 
that is $\hat\Sg_s^+$ consists of 
infinite sequences $\om=\ld\om_3\om_2\om_1$ of elements from the set
$\{1,2,\ld,s\}$. 
For each $\om \in \hat\Sigma _{s}^{+}$, 
by $\om|_n$ we
denote the finite word $(\om_n,\om_{n-1},\ld,\om_2,\om_1)$; more generally,
for $b\ge a$ we put $\om|_b^a=(\om_b,\om_{b-1},\ld,\om_a)$. 
For each finite word $\om=(\om_{n},\ldots\!,\om_{1})$,\ we set 
$f_{\om}:=f_{\om_{1}}\circ \cdots \circ f_{\om_{n}}.$ 
The technical 
tool that allows us to develop the further machinery is the following.

\

\blem\lab{l4042805}
Let $G=\langle f_1,\ld,f_s\rangle $ be a rational semigroup satisfying the condition 
(E1). Let $R_{0}$ be the number in Lemma~\ref{l3042805}. 
Fix an integer $q\ge 1$and a real number $\l\in(0,1)$. Then for every finite
set $E\sbt J(G)$, every $\om\in \hat\Sg_s^+$, every \\ 
$R\in \\ 
\lt(0,\min\lt\{1,\frac{1}{2}R_{0}, \ka^{-1/2}, \frac{1}{2}\dist\(E,CV(f_{\om|_q})\),
4^{-1}\min\{d(z,\xi):z,\xi\in E, z\ne\xi\}\rt\}\rt),
$\\  
every integer $n\ge 0$,
and every $z\in E$, there exists a subset $I_n(z,\om)$ of the set of 
all inverse meromorphic 
branches of $f_{\om|_{qn}}$ defined on $B(z,2R)$ and satisfying the following
properties with $I_n=\bu_{z\in E}I_n(z,\om)$.
\begin{itemize}
\item[($a_n$)] If $z\in E$ and $\phi\in I_{n+1}(z,\om)$, then 
$f_{\om|_{q(n+1)}^{qn+1}}\circ\phi\in I_n(z,\om)$.
\item[($b_n$)] If $\phi\in I_n$, then $\diam(\phi(B(z,R)))
\le K_{G}(2R)^2\ka^{-1/2}\l^{n/2}C_{sa}^{1/2}$.
\item[($c_n$)] $\phi(B(z,2R))\cap CV\(f_{\om|_{q(n+1)}^{qn+1}}\)=\es$ for
all $z\in E$ and all $\phi\in I_n$.
\item[($d_n$)] $\#(J_n\sms I_n)\le dq+\l^{-n}$ for all $n\ge 1$, where $J_n$
is the family of all compositions of all maps $\phi\in I_{n-1}(z,\om)$,
$z\in E$, with all meromorphic inverse branches of $f_{\om|_{qn}^{q(n-1)+1}}$.
\item[($e$)] $I_{0}=\{ \mbox{{\em Id}}|_{B(z,2R)}\mid z\in E\} .$  
\end{itemize}
\elem
{\sl Proof.} We shall construct recursively the sets $I_n(z,\om)$, $n\ge 0$,
such that the conditions ($a_n$), ($b_n'$), ($c_n$) and ($d_n$) 
(here $n\ge 1$), where ($b_n'$) requires that 
\begin{itemize}
\item[($b_n'$)] If $\phi\in I_n$, then $A\(\phi(B(z,R))\)\le \l^n$. 
\end{itemize}
The base of induction, the family $I_0$ consists of all the identity maps
defined on the balls $B(z,2R)$, $z\in E$. The condition ($b_0'$) is satisfied
since $A(\oc)=1$ and ($c_0$) is satisfied because of the choice of the radius
$R$. Now assume that for some $n\ge 0$ the subsets $I_n(z,\om)$, $z\in E$,
have been constructed so that the conditions ($b_n'$) and ($c_n$) are satisfied.
The inductive step is to construct the subsets $I_{n+1}(z,\om)$, $z\in E$,
so that the conditions ($a_{n+1}$), ($b_{n+1}'$), ($c_{n+1}$) and ($d_{n+1}$) are satisfied.
This will complete our recursive construction. In view of ($c_n$) all the
meromorphic inverse branches of $f_{\om|_{q(n+1)}^{qn+1}}$ are well defined on all
the sets $\phi(B(z,2R))$, $z\in E$, $\phi\in I_n(z,\om)$. Their compositions 
with corresponding elements $\phi\in I_n(z,\om)$ are said to form the subset 
$J_{n+1}(z,\om)$. Note that $J_{n+1}=\bu_{z\in E}J_{n+1}(z,\om)$. The subset 
$I_{n+1}(z,\om)$, $z\in E$, is defined to consist of all the elements $\psi
\in J_{n+1}(z,\om)$ for which the following two conditions are satisfied.
\begin{itemize}
\item[(i)] $A(\psi(B(z,R)))\le \l^{n+1}$.
\item[(ii)] $\psi(B(z,2R))\cap CV\(f_{\om|_{q(n+2)}^{q(n+1)+1}}\) =\es$.
\end{itemize}
Thus, conditions ($b_{n+1}'$) and ($c_{n+1}$) are satisfied immediately. 
Condition ($a_n$) is satisfied since it holds for all $\psi\in 
J_{n+1}(z,\om)$, $z\in E$, and $I_{n+1}(z,\om)$ is a subset of $J_{n+1}(z,\om)$.
We are left to show that ($d_{n+1}$) holds. Now, if $\psi\in J_{n+1}$, 
say $\psi\in J_{n+1}(z,\om)$, $z\in E$, but (i) is not satisfied, then
$A(\psi(B(z,R)))> \l^{n+1}$. Since all the sets $\psi(B(z,R))$, 
$\psi\in J_{n+1}(z,\om)$, $z\in E$, are mutually disjoint and since $A(\oc)=1$,
we conclude that the number of elements of $J_{n+1}$ for which 
condition (i) fails is bounded above by $1/ \l^{n+1}=\l^{-(n+1)}$. 
Since the number of critical points of each generator of the semigroup
$G$ is bounded above by $d$, the cardinality of the set of critical values 
of $f_{\om|_{q(n+2}^{q(n+1)+1}}$ is bounded above by $dq$. Since all the 
sets $\psi(B(z,2R))$, $\psi\in J_{n+1}(z,\om)$, $z\in E$, are mutually 
disjoint, we thus conclude that the number of elements $\psi\in J_{n+1}$
for which condition (ii) fails, is bounded above by $dq$. In conclusion
$\#(J_{n+1}\sms I_{n+1})\le dq+\l^{-(n+1)}$, meaning that ($d_{n+1}$) is 
satisfied. The recursive construction is complete. Since 
$G$ satisfies (E1), it follows from Lemma~\ref{l3042805} that
for all $n\ge 0$, all $z\in E$, and all $\phi\in I_n(z,\om)$, we have
$$
B\(\phi(z),K_{G}(2R)^{-1}|\phi'(z)|R\)
\sbt \phi(B(z,R))
\sbt 
B\(\phi(z),K_{G}(2R)|\phi'(z)|R\).
$$
Hence, making use of ($b_n'$), we get that
$$
\aligned
\diam^2\(\phi(B(z,R))\)
&\le \(2K_{G}(2R)|\phi'(z)|R\)^2\\ 
& =K_{G}(2R)^4\ka^{-1}\(2\ka ^{1/2}K_{G}(2R)^{-1}|\phi'(z)|R\)^2\\
&=K_{G}(2R)^4\ka^{-1}C_{sa}A\(B\(\phi(z),K_{G}(2R)^{-1}|\phi'(z)|R\)\)\\ 
& \le K_{G}(2R)^4\ka^{-1}C_{sa}A\(\phi(B(z,R))\) \\
&\le K_{G}(2R)^4\ka^{-1}C_{sa}\l^n.
\endaligned
$$
Thus, $\diam\(\phi(B(z,R))\)\le K_{G}(2R)^2\ka^{-1/2}C_{sa}^{1/2}\l^{n/2}$.
We are done. \endpf

\

\ni In order to simplify the notation, put 
$$
\aligned
\ & R(E,\om|_q):=\\ 
&\min\lt\{1,\frac{1}{2}R_{0}, \ka^{-1/2},\frac{1}{2}\dist\(E,CV\(f_{\om|_q}\)\),
             4^{-1}\min\{d(z,\xi):z,\xi\in E, z\ne\xi\}\rt\}.
\endaligned
$$

\

\bcor\lab{c5042905}
Suppose that $G=\langle f_1,\ld,f_s\rangle $ is an E-semigroup of rational maps. 
Fix an integer $q\ge 1$and a real number $\l\in(0,1)$.   
Then for every $z\in J(G)$,\  there exists a number $R=R_q(z)>0$, 
a number $R'_{q}(z)>0$, 
and a number $D\geq 1$, where $D$ does not depend on 
$q,\l ,z$,\     
such that for
every $\om\in \hat\Sg_s^+$ and every integer $n\ge 1$, there exists
$W_n(z,\om),\ Z_{n}(z,\om )$, a subset of the set of 
all connected components of $f_{\om|_{qn}}^{-1}(B(z,R))$,
with the following properties.
\begin{itemize}
\item[($A_n$)] If $V\in W_{n+1}(z,\om)$, then $f_{\om|_{q(n+1)}^{qn+1}}(V)\in
W_n(z,\om)$.
\item[($B_n$)] If $V\in W_{n}(z,\om)$, then $\diam(V)
\le K_{G}(R'_{q}(z))^4\ka^{-1/2}\l^{n/2}C_{sa}^{1/2}$.
\item[($C_n$)] $Z_n(z,\om)$ is the family of all connected components of 
the sets $f_{\om|_{qn}^{q(n-1)+1}}^{-1}(V)$, $V\in W_{n-1}(z,\om)$, and  
$\#\(Z_n(z,\om)\sms W_n(z,\om)\)\le dq+\l^{-n}$, and 
\item[($D_n$)] For every $V\in Z_n(z,\om)$, the map $f_{\om|_{qn}}|_V:V\to\oc$
is at most $D$-to-$1$.
\item[($E_{1}$)] $W_{1}(z,\om )=Z_{1}(z,\om )=
\{ \mbox{connected components } V \mbox{ of }f_{\om |_{q}}^{-1}(B(z,R))\} .$ 
\end{itemize}
Note also that in fact $W_n(z,\om)$ depends only on $z$ and $\om|_{qn}$, so we
can and will in the forthcoming sections write $W_n(z,\om|_{qn})$
\ecor
{\sl Proof.}  
Let $z\in J(G)$ be a point. Since 
$J(G)=\pi _{2}(J(f))$,\ there exists a point 
$\om \in \Sigma _{s}$ such that 
$(\om ,z)\in J(f).$ Then, 
from the assumption E2,\ E3,\ 
there 
exists $p\ge 1$ independent of $q$ (but depending on $(\om ,z)$) such that
$$ 
\( \bigcup _{r\in \Bbb{N}}(f^{pq+r})^{-1}(\om ,z) \) \cap 
\mbox{\Crit}(f)=\emptyset .
$$ 
Then, we obtain 
$$
\bu_{\rho\in\Sg_s^*}\bu_{|\tau|=pq}f_\rho^{-1}(f_\tau^{-1}(z))\cap 
\Crit(f_{\rho })=\es.
$$
In particular
\begin{equation}\label{1061505}
\lt(\bu_{|\tau|=pq}f_\tau^{-1}(z)\rt)\cap CV(f_\rho)=\es
\end{equation}
for all $\rho\in\Sg_s^*$. Set $E=\bu_{|\tau|=pq}f_\tau^{-1}(z)$. It follows
from (\ref{1061505}) that 
$$ \hat R_q(z)=\frac{1}{2}\lambda ^{p/2}\cdot 
\min\{R(E,\rho):\rho\in\{1,2,\ld,s\}^{q}\} >0. $$

Now, there is $R_q(z)>0$ so small that the following two conditions are satisfied.
\begin{itemize}
\item[(a)] 
For each $\tau\in\{1,2,\ld,s\}^{pq}$ each connected component of 
$f_\tau^{-1}(B(z,R_q(z)))$ is contained in exactly one ball
$B(\xi,\hat R_q(z))$, where $\xi\in E$.
\item[(b)] For each $\g\in \Sg_s^*$ with $|\g|\le pq$, each connected component of 
$f_\g^{-1}(B(z,R_q(z)))$ has the diameter bounded above by 
$\ka^{-1/2}\l^{|\g|/2}$.
\end{itemize}
Now, for every $1\le k\le p$ and every $\om\in \hat \Sg_s^+$, define 
$W_k(z,\om)$ and $Z_k(z,\om )$ to be
the family of all connected components of
$f_{\om|_{qk}}^{-1}(B(z,R_q(z)))$. The conditions 
($A_k$), ($B_k$), and ($C_k$) are obviously satisfied for all $1\le k\le p$. Now, for
every $\om\in \hat \Sg_s^+$, every $n\ge p+1$ and every $\xi\in
f_{\om|_{pq}}^{-1}(z)\sbt E$, 
consider the inverse branch $\phi:B(\xi,2\lambda ^{-p/2}\hat R_q(z))\to\oc\in 
I_{n-p}(\xi,\om|_\infty^{pq})$.
It follows from Lemma~\ref{l4042805}($b_{n-p}$) and Lemma~\ref{l3042805} that
$$
\aligned
2K_{G}(2\lambda ^{-\frac{p}{2}}\hat{R}_{q}(z))^{-1} \lambda
^{-p/2}|\phi'(\xi)|\hat R_q(z) 
&\le\diam\(\phi(B(\xi,\lambda ^{-p/2}\hat R_q(z)))\) \\
&\le K_{G}(2\lambda ^{-\frac{p}{2}}\hat{R}_{q}(z))^{2}\ka^{-1/2}\l^{{n-p\over 2}}C_{sa}^{1/2}.
\endaligned
$$
Thus $|\phi'(\xi)|\le 
\frac{1}{2}K_{G}(2\lambda ^{-\frac{p}{2}}\hat{R}_{q}(z))^{3}
\ka ^{-1/2}\lambda ^{n/2}\hat R_{q}(z)^{-1}C_{sa}^{1/2}$
, and therefore,
using Lemma~\ref{l3042805} again and the definition of $\hat R_q(z)$, we get that 
$$
\diam\(\phi(B(\xi,\hat R_q(z)))\)
\le K_{G}(2\lambda ^{-\frac{p}{2}}\hat{R}_{q}(z))^{4}\ka^{-1/2}\l^{n/2}C_{sa}^{1/2}.
$$
So, looking also at Lemma~\ref{l4042805}($d_n$) and ($c_n$), we complete the proof 
of items ($A_n$), ($B_n$) and ($C_n$) by defining for every $n\ge p+1$ the family 
$W_n(z,\om)$ to consist of all the sets of the form $\phi(V_\xi)$, where $\xi\in
f_{\om|_{pq}}^{-1}(z)$, $V_\xi$ is a connected component of 
$f_{\om|_{pq}}^{-1}(B(z,R_q(z)))$ contained in $B(\xi,\hat R_q(z))$, and $\phi\in
I_{n-p}(\xi,\om_\infty^{pq})$. Decreasing $R_q(z)$ appropriately, the items 
($D_n$) follow now from this construction and  Lemma~\ref{l4042805}. 
\endpf

\

\section{Perron-Frobenius Operators and Gibbs States}

\ni From now on throughout the entire paper assume that $\psi:J(f)\to\R$
is a H\"older continuous function. Given $n\ge 1$, $\om\in\{1,\ld,s\}^n$,
and a continuous function $g:J(f)\to\R$, define
$$
\L_{\psi,\om}^{(n)}g(\tau,z)
=\sum_{x\in f_\om^{-1}(z)}\exp\(S_n\psi(\om\tau,x)\)g(\om\tau,x),
$$
where here and in the sequel the summation is taken with multiplicities
of all critical points of $f_\om$, and 
$S_{n}\psi := \sum _{j=0}^{n-1}\psi \circ f^{j}.$ 
 Since $\sup\{\#(f_\om^{-1}(z)):z\in \C\}
<\infty$, $\L_{\psi,\om}^{(n)}$ is a bounded linear operator acting on the
Banach space $C(J(f))$ of continuous functions on $J(f)$ endowed with the
supremum norm. Set
\begin{equation}
\label{pfopdf}
\L_\psi^{(n)}:=\sum_{|\om|=n}\L_{\psi,\om}^{(n)},\  
\L_\psi :=\L_{\psi }^{(1)}.
\end{equation}
Then $\L_\psi^{(n)}$ and $\L_{\psi }$ also 
act continuously on the Banach space $C(J(f)).$ 
We call $\L_{\psi }:C(J(f))\rightarrow C(J(f))$ the 
{\bf Perron-Frobenius operator associated with the potential 
$\psi .$}  
Note that 
$$
\L_\psi^{(n)}g(x)
=\sum_{y\in f^{-n}(x)}\exp\(S_n\psi(y)\)g(y)
=\L_\psi^ng(x).
$$
Define the pointwise pressure $\P ^{p}(\psi)$
of the function $\psi$ by the following formula.
\begin{equation}
\lab{PxPp}
\aligned
\P_x(\psi)
:=\limsup_{n\to\infty}{1\over n}\log\L_\psi^n\1(x), \  x\in J(f), \ 
\text{ and } \\  
\P^p(\psi):=\sup\{\P_{x}(\psi):x\in J(f)\} , \ \ \ \ \ \ \ \ \ \ \ \ \ \ \ \ \ \ \ \  
\endaligned
\end{equation}
where $\1(x) :=1.$ 
Throughout the entire paper we work with the assumption that
\begin{equation}\lab{6061805}
\P^p(\psi)>\sup(\psi)+\log s.
\end{equation}
In particular, we can fix a point $b\in J(f)$ such that
\begin{equation}\lab{2061605}
\eta:=\exp\(\sup(\psi)+\log s-\P_b(\psi)\)<1.
\end{equation}
Fix $\l\in(0,1)$. There then exists $q=q(\l)\ge 1$ such that
\begin{equation}\lab{1061705}
\g_q^{-1}:=dq+\l^{-1}\le (1+2D)^{-1}\eta^{-q}.
\end{equation}
Let $\L_\psi^*:C(J(f))^*\to C(J(f))^*$ be the operator conjugate to $\L_\psi$,
i.e., $\L_\psi^*\nu(g)=\nu(\L_\psi g)$.
Using (\ref{2061605}) and the fact that the map $f:J(f)\to J(f)$ is open,
as an immediate consequence of Theorem~3.9 and Proposition 2.2 
in \cite{duecm}, we get the following.

\

\blem\lab{l1061605}
If $G=\langle f_1,\ld,f_s\rangle $ is an E-semigroup of rational maps, $\psi:J(f)\to
\R$ is a H\"older continuous potential satisfying (\ref{6061805}) and $b\in J(f)$
is selected so that (\ref{2061605}) holds, then there exists a Borel atomless
probability measure $m_\psi$ on $J(f)$ such that $\L_\psi^*m_\psi
=e^{\P_b(\psi)}m_\psi$.
\elem

\

\ni The measure $m_\psi$ is called $\exp(\P_b(\psi)-\psi)$-conformal for $f$.
Obviously $(\L_\psi^*)^nm_\psi=e^{\P_b(\psi)n}m_\psi$ for all $n\ge 0$, and this 
equivalently means that
$$
m_\psi(f^n(A))=\int_A\exp\(\P_b(\psi)n-S_n\psi\)dm_\psi
$$ 
for every Borel set $A\sbt J(f)$ for which the restriction $f^n|_A$ is injective.

\  

\brem\lab{rmeas}
{\rm Note that all forthcoming considerations depend only on the 
above relation and not on the particular way the measure 
$m_\psi $ was constructed.  }
\erem 

\ 

From now on throughout the paper put
$$
\ov\psi=\psi-\P_b(\psi).
$$
Now set
$$
L_\om^{(n)}=\L_{\ov\psi,\om}^{(n)}, \  \L^n=\L_{\ov\psi}^n,  \
\text{ and } \  \L=\L_{\ov\psi}.
$$
Now for every $z\in J(G)$, every $n\ge 1$, every $\om\in\{1,\ld,s\}^{qn}$, and
every $g\in C(J(f))$, define the function 
$G_{z,\om}^{n}g:\pi_2^{-1}(B(z,R_q(z)))
\to\R$ by setting
\begin{equation}\lab{1061605}
G_{z,\om}^ng(\tau,\xi)
=\sum_{V\in W_n(z,\om)}\sum_{x\in f_\om^{-1}(\xi)\cap V}
 \exp\(S_{qn}\ov\psi(\om\tau,x)\)g(\om\tau,x),
\end{equation}
where $R_q(z)$ and $W_{n}(z,w)$ come from Corollary~\ref{c5042905}. Since for every $V\in
W_n(z,\om)$ the map $f_{\om }:V\to B(z, R_q(z))$ is a branched covering,
for every $\xi\in B(z,R_q(z))$ there is a bijection 
$\hat\xi:f_{\om }^{-1}(z)
\cap V\to f_{\om }^{-1}(\xi)\cap V$, where all critical points of 
$f_{\om }$
in $f_{\om }^{-1}(z)$ and $f_{\om }^{-1}(\xi)$ are counted with multiplicities.
So, (\ref{1061605}) can be rewritten in the following form.
$$
G_{z,\om}^ng(\tau,\xi)
=\sum_{V\in W_n(z,\om)}\sum_{x\in f_{\om }^{-1}(z)\cap V}
 \exp\(S_{qn}\ov\psi(\om\tau,\hat\xi(x))\)g(\om\tau,\hat\xi(x)).
$$
Since the function $\psi:J(f)\to\R$ is H\"older continuous, it follows from
Corollary~\ref{c5042905}($B_n$) that there exists a 
constant $H>0$
such that (with $x\in f_{\om }^{-1}(z)\cap V$),
$$
|S_{qn}\ov\psi(\om\tau,\hat\xi(x))-S_{qn}\ov\psi(\om\th,x)|\le H
$$
for all $\tau,\th\in\Sg_s$, or equivalently,
\begin{equation}\lab{6081605}
e^{-H} \exp\(S_{qn}\ov\psi(\om\th,x)\)
\le \exp\(S_{qn}\ov\psi(\om\tau,\hat\xi(x))\)
\le e^{H} \exp\(S_{qn}\ov\psi(\om\th,x)\).
\end{equation}
In consequence
\begin{equation}\lab{2061605a}
e^{-H} G_{z,\om}^n\1(\th,z)
\le G_{z,\om}^n\1(\tau,\xi)
\le e^{H} G_{z,\om}^n\1(\th,z).
\end{equation}
Set 
$$
G_z^n=\sum_{|\om|=qn}G_{z,\om}^n.
$$
It then follows from (\ref{2061605a}) that
\begin{equation}\lab{4061605}
e^{-H} G_z^n\1(\th,z)
\le \sum_{|\om|=qn}G_{z,\om}^n\1 (\tau _\om ,\xi _\om )
\le e^{H} G_z^n\1(\th,z)
\end{equation}
for all $n\ge 0$, all $(\th,z)\in J(f)$ and all $(\tau_\om,\xi_\om)\in 
\pi_2^{-1}(B(z,R_q(z)))$. 
Since $J(G)$ is a compact set, there exist finitely many points, say 
$z_1,z_2,\ld,z_u\in J(G)$ such that $\bu_{j=1}^uB\(z_j,R_q(z_j))\spt J(G)$.
Put 
$$
M_{j,\om }^{n}g=L_{\om }^{(qn)}|_{\pi_2^{-1}\(B\(z_j,R_q(z_j)\)\)}, \ 
M_{j}^{n}g=\sum _{|\om |=qn}M_{j,\om }^{n}g, \ 
G_{j,\om}^n=G_{z_j,\om}^n, \
G_j^n=G_{z_j}^n,
$$
and 
$$
\aligned
|||G_j^{n}g|||_\infty =\sum\limits_{|\om|=qn}||G_{j,\om}^{n}g||_\infty, \
|||M_j^{n}g|||_\infty =\sum\limits_{|\om|=qn}||M_{j,\om}^{n}g||_\infty, \\ 
|||\L^{n}g|||_\infty =\sum\limits_{|\om|=n}||L_{\om }^{(n)}g||_\infty.\ \ \ \ \ 
\ \ \ \ \ \ \ \ \ \ \ \ \ \ \ \ \ \ \ \ \ \ \ \ \ \ \ \ \ \ \ \ \ \ \ \ 
\endaligned
$$ 
We shall prove the following. 

\ 

\blem\lab{l2061605} It holds
$$
0<Q_\psi:=\max_{1\le j\le u}\sup_{n\ge 0}\{|||G_j^n\1|||_\infty\}<\infty.
$$
\elem
{\sl Proof.} Since the map $f:J(f)\to J(f)$ is topologically exact and since
$\{B\(z_j,R_q(z_j)\)\}_{j=1}^u$ is an open cover of $J(f)$, there exists $k\ge 1$
such that for all $j=1,2,\ld,u$,
\begin{equation}\lab{3061605}
f^{kq}\(\pi_2^{-1}\(B\(z_j,R_q(z_j)\)\)\)\spt J(f).
\end{equation}
By Lemma~\ref{l1061605}, $\L^{*n}m_\psi=m_\psi$ for all $n\ge 0$, and
consequently $\int\L^n\1 dm_\psi=\int\1 dm_\psi=1$. Fix $n\ge 0$. There then 
exists $w_0\in J(f)$ such that
\begin{equation}\lab{5061605}
\L^{q(k+n)}\1(w_0)\le 1.
\end{equation}
Now fix an arbitrary $1\le j\le u$, an arbitrary $\om\in\{1,\ldots\! ,s\}^{qn}$,
and an arbitrary $x_\om\in \pi_2^{-1}\(B\(z_j,R_q(z_j)\)\)$. By 
(\ref{3061605}) there exists $y_j\in \pi_2^{-1}\(B\(z_j,R_q(z_j)\)\)$
such that $f^{kq}(y_j)=w_{0}$. Applying (\ref{4061605}) with $z=z_j$, we get that
$$
\sum_{|\om|=qn}G_{j,\om}^n\1(x_\om) \le e^{2H}G_j^n\1(y_j).
$$
Also, by (\ref{5061605}), we obtain
$$
\aligned
G_j^n\1(y_j)
& \le \L^{qn}\1(y_j)
\le ||\exp\(S_{kq}(-\ov\psi)\)||_\infty\L^{q(k+n)}\1(w_{0})\\ 
&  \le ||\exp\(S_{kq}(-\ov\psi)\)||_\infty.
\endaligned
$$
Thus, $\sum_{|\om|=qn}G_{j,\om}^n\1(x_\om) \le e^{2H}||
\exp\(S_{kq}(-\ov\psi)\)||_\infty$, and we are done by taking supremum over 
all $x_\om\in \pi_2^{-1}\(B\(z_j,R_q(z_j)\)\)$. \endpf

\

\ni Now we are in position to prove the required upper bound on the 
iterates of the Perron-Frobenius operator $\L$.

\

\blem\lab{l1061705}
For every $n\ge 1$, we have that 
$$
\aligned
|||\L^{qn}\1|||_\infty
& \le Q_\psi+D\sum_{k=1}^n\(\eta^q\g_q^{-1}\) ^k|||\L^{q(n-k)}\1|||_\infty \\ 
& \le Q_\psi+D\sum_{k=1}^n\(1+2D\)^{-k}|||\L^{q(n-k)}\1|||_\infty.
\endaligned
$$
\elem
{\sl Proof.} Fix $n\ge 1$, 
$\om=(\om_{qn},\om_{qn-1},\ld,\om_1)\in\{1,\ld,s\}^{qn}$,
and $(\tau,x)\in J(f)$. 
There then exists $1\le j\le u$ such that 
$(\tau,x)\in \pi_2^{-1}\(B\(z_j,R_q(z_j)\)\)$. 
For each $a,b \in \Bbb{N}$ with $a\leq b\leq qn$, 
we set $\om |_{b}^{a}:=(\om _{b},\ldots\! ,\om _{a}).$ Moreover, 
for each $l\in \Bbb{N}$ with $l\leq qn$, 
we set $\om |_{l}=(\om _{l},\ldots\! ,\om _{1}).$ 

One can now represent
$L_\om^{(qn)}\1(\tau,x)$ in the following way.
\begin{equation}\lab{3081605}
\aligned
\ & \ \ \ \ M_{j,\om}^{n}\1(\tau,x)\\ 
&=   \sum_{k=1}^n\sum_{V\in Z_{k}(z_j,\om)\sms W_{k}(z_j,\om)}
     \sum_{y\in V\cap f_{\om|_{qk}}^{-1}(x)}\exp\(S_{qk}\ov\psi(\om|_{qk}\tau,y)\)
     L_{\om|_{qn}^{qk+1}}^{(q(n-k))}\1(\om|_{qk}\tau,y) \\
&{}\  \  \  \  \  \  \  \  \     + G_{j,\om}^{n}\1(\tau,x) \\
&\le \sum_{k=1}^n\sum_{V\in Z_{k}(z_j,\om)\sms W_{k}(z_j,\om)}
     \sum_{y\in V\cap f_{\om|_{qk}}^{-1}(x)}\exp\(qk\sup(\ov\psi)\)
     L_{\om|_{qn}^{qk+1}}^{(q(n-k))}\1(\om|_{qk}\tau,y) \\ 
&{}\ \ \ \ \ \ \ \ \ + ||G_{j,\om}^n\1||_\infty.
\endaligned
\end{equation}
Applying now Corollary~\ref{c5042905} ($C_{k}$) and ($D_{k}$), we 
estimate further as follows.
$$
M_{j,\om}^{n}\1(\tau,x)
\le \sum_{k=1}^nD(dq+\l^{-k})\exp\(qk\sup(\ov\psi)\)
\| L_{\om|_{qn}^{qk+1}}^{(q(n-k))}\1\| _{\infty } + ||G_{j,\om}^n\1||_\infty.
$$
Since $dq+\l^{-k}\le \g_q^{-k}$ (see (\ref{1061705})), we thus get
$$
M_{j,\om}^{n}\1(\tau,x)
\le D\sum_{k=1}^n\(\exp\(q\sup(\ov\psi)\)\g_q^{-1}\)^k
 ||L_{\om|_{qn}^{qk+1}}^{(q(n-k))}\1||_\infty + ||G_{j,\om}^n\1||_\infty.
$$
Taking supremum over all $(\tau,x)\in \pi_2^{-1}\(B\(z_j,R_q(z_j)\)\)$, we
thus get
$$
||M_{j,\om}^{n}\1||_\infty
\le D\sum_{k=1}^n\(\exp\(q\sup(\ov\psi)\)\g_q^{-1}\)^k
    ||L_{\om|_{qn}^{qk+1}}^{(q(n-k))}\1||_\infty + ||G_{j,\om}^n\1||_\infty.
$$
So, summing over all words $\in\{1,\ld,s\}^{qn}$, we obtain using 
Lemma~\ref{l2061605}, the following.
\begin{equation}\lab{3061805}
\aligned
|||M_j^{n}\1|||_\infty
&\le D\sum_{k=1}^n\(\exp\(q\sup(\ov\psi)\)\g_q^{-1}\)^ks^{qk}
     |||\L^{q(n-k)}\1|||_\infty + |||G_j^n\1|||_\infty \\
&\le D\sum_{k=1}^n\exp\(qk(\sup(\ov\psi)+\log s)\)(\g_q^{-1})^k
     |||\L^{q(n-k)}\1|||_\infty + Q_\psi \\
&=   D\sum_{k=1}^n\(\eta^q\g_q^{-1}\)^k|||\L^{q(n-k)}\1|||_\infty + Q_\psi.
\endaligned
\end{equation}
Since this inequality holds for all $j=1,\ld,u$, we thus get
$$
|||\L^{qn}\1|||_\infty
\le D\sum_{k=1}^n\(\eta^q\g_q^{-1}\)^k|||\L^{q(n-k)}\1|||_\infty + Q_\psi.
$$
The first inequality to be proved is thus established. In order to derive 
the second one from it, invoke (\ref{1061705}). \endpf

\

\blem\lab{l1061805}
There exists a constant $\ov Q_\psi>0$ such that for all $n\ge 0$, we have
$$
|||\L^n\1|||_\infty\le \ov Q_\psi.
$$
\elem
{\sl Proof.} We first shall prove by induction that
\begin{equation}\lab{1061805}
|||\L^{qn}\1|||_\infty\le 2Q_\psi
\end{equation}
for all integers $n\ge 0$. Since $||\1||_\infty=1$, this formula holds for
$n=0$. So, fix $n\ge 1$ and suppose that (\ref{1061805}) is true for all 
$0\le k\le n-1$. It then follows from Lemma~\ref{l1061705} that
$$
\aligned
|||\L^{qn}\1|||_\infty\
& \le Q_\psi+D\sum_{k=1}^n\(1+2D\)^{-k}2Q_\psi
=Q_\psi\lt(1+2D\sum_{k=1}^n\(1+2D\)^{-k}\rt)\\ 
& \le 2Q_\psi,
\endaligned
$$
and (\ref{1061805}) is proved. 
Since $|||\L^{i+j}\1|||_\infty
\le 
|||\L^i\1|||_\infty |||\L^j\1|||_\infty$, we are
done by setting $\ov Q_\psi=2Q_\psi\max\{|||\L^k\1|||_\infty:k=0,1,\ld,q-1\}$.
\endpf

\

\blem\lab{l1061805a}
There exists a constant $\un Q_\psi>0$ such that for all $n\ge 0$, we have
$$
\inf _{y\in J(f)}\L^n\1 (y)\ge \un Q_\psi.
$$
\elem
{\sl Proof.} Taking $q\ge 1$ sufficiently large, we can make the product
$\eta^q\g_q^{-1}$ (see (\ref{1061705})) as small as we wish. It therefore 
follows from (\ref{3061805}) and Lemma~\ref{l1061805} that for $q\ge 1$
large enough, for all $n\ge 0$
and all
$1\le j\le u$, we have
\begin{equation}\lab{4061805}
|||M_j^{n}\1|||_\infty\le {1\over 2}+|||G_j^n\1|||_\infty.
\end{equation}
Since, by Lemma~\ref{l1061605}, $\int\L^{qn}\1dm_\psi=\int\1 dm_\psi=1$, 
there thus exists $x\in J(f)$ such that $\L^{qn}\1(x)\ge 1$. Since 
$\{\pi_2^{-1}(B(z_j,R_q(z_j)))\}_{j=1}^u$ is a cover of $J(f)$, there exists
$1\le i\le u$ such that $x\in \pi_2^{-1}(B(z_i,R_q(z_i)))$. It follows from
(\ref{4061805}) that $|||G_i^n(\1)|||_\infty\ge 1/2$. Applying now 
(\ref{4061605}) we see that
\begin{equation}\lab{5061805}
G_i^n\1(w)\ge (2e^{2H})^{-1}
\end{equation}
for all $w\in \pi_2^{-1}(B(z_i,R_q(z_i)))$. Take now an arbitrary point $y\in J(f)$.
With $k\ge 1$, as in the proof of Lemma~\ref{l2061605}, it follows from 
(\ref{3061605}) that there exists $\ov y\in\pi_2^{-1}(B(z_i,R_q(z_i)))$ such
that $f^{kq}(\ov y)=y$. So, using (\ref{5061805}) and the definition of the
Perron-Frobenius operator, we obtain
$$
\aligned
\L^{q(n+k)}\1(y)
&\ge \exp\(S_{kq}(\ov\psi(\ov y))\)\L^{qn}\1(\ov y)
 \ge \exp\(\inf\(S_{kq}(\ov\psi)\)\)G_i^n\1(\ov y) \\
&\ge (2e^{2H})^{-1}\exp\(\inf\(S_{kq}(\ov\psi)\)\)
\endaligned
$$
Since $\inf\(\L^{a+b}\1\)\ge \inf\(\L^a\1\)\inf\(\L^b\1\)$, we are therefore 
done by taking 
$$
\aligned
\un Q_\psi:=\ \ \ \ \ \ \ \ \ \ \ \ \ \ \ \ \ \ \ \ \ \ \ \ \ \ \ \ \ \ 
\ \ \ \ \ \ \ \ \ \ \ \ \ \ \ \ \ \ \ \ \ \ \ \ \ \ \ \ \ \ \ \ \ \ 
\ \ \ \ \ \ \ \ \ \ \ \ \ \ \ \ \ \ \ \ \ \ \ \ \ \\ 
\min \left\{ (2e^{2H})^{-1}\exp\(\inf\(S_{kq}(\ov\psi)\)\) , 1\right\} 
\cdot 
\min\left\{\inf\(\L^u\1\):u=0,1,\ld,kq-1\right\}.
\endaligned
$$
\endpf

\

\ni Now repeating verbatim the proof of Lemma~20 from \cite{dues}, using
Lemma~\ref{l1061805} and 
Lemma~\ref{l1061805a}, we get the following.

\

\blem\lab{l3061805} 
Let $\tilde{h}(x):=\liminf _{n\rightarrow \infty }
\frac{1}{n}\L ^{n}\1 (x)$, 
$h_{1}(x):=\inf _{n\geq 0}\L ^{n}\tilde{h}(x)$ and 
$h(x):= \frac{h_{1}(x)}{\int _{J(f)}h_{1}\ dm_{\psi }}$ for 
each $x\in J(f).$ 
Then 
$h:J(f)\rightarrow \Bbb{R} $ is a Borel measurable function 
such that the following hold.
\begin{itemize}
\item[(a)]
$$
{\un Q_\psi\over \ov Q_\psi}\le h(x)\le{\ov Q_\psi\over \un Q_\psi} \  
\text{ for every } \  x\in J(f).
$$
\item[(b)] $\L h(x)=h(x)$ for every $x\in J(f)$.
\item[(c)] $\int hdm_\psi=1$.
\end{itemize}
\elem

\

\brem\lab{rnormal}   
Note that up to a normalized factor ( to make (c) hold ) 
the function $h$ is independent of the conformal measure 
$m_{\psi }.$  
\erem

\ 

\ni As an immediate consequence of this lemma and Proposition 2.2 
in \cite{duecm}, we have the following.

\

\bthm\lab{t4061805}
The Borel probability measure $\mu_\psi=hm_\psi$ is $f$-invariant, equivalent 
to $m_\psi$, and ${d\mu_\psi\over dm_\psi}
\in[\un Q_\psi/\ov Q_\psi, \ov Q_\psi/ \un Q_\psi]$.
\ethm

\section{Equilibrium States}

\ni Our objective in this section is to show that the measure $\mu_\psi$ 
produced in Theorem~\ref{t4061805} is a unique equilibrium state for the
potential $\psi$ and that it is ergodic. Let $\P(\psi)$ be the ordinary
topological pressure of the potential $\psi$. If $\mu$ is a Borel
probability $f$-invariant measure on $J(f)$, denote by $J_\mu$ its Jacobian
with respect to the map $f$ (see page 108 in \cite{P}). We start with the following.

\

\blem\lab{l5061805}
$\P(\psi)\ge \P_b(\psi)$.
\elem
{\sl Proof.} Since $h_{\mu_\psi}(f)\ge \int_{J(f)}\log J_{\mu_\psi}d\mu_\psi$
(this is true for every finite-to-one endomorphisms and every probability invariant measure. 
See Lemma 10.5 and Theorem 5.14 in \cite{P})
and since $J_{\mu_\psi}={h\circ f\over h}\exp(\P_b(\psi)-\psi)$ everywhere, 
it follows from Theorem~\ref{t4061805} and the Variational Principle that 
\begin{equation}\lab{2081605}
\P(\psi)
\ge h_{\mu_\psi}(f)+\int\psi d\mu_\psi
\ge \int(\P_b(\psi)-\psi+\psi)d\mu_\psi
=\P_b(\psi).
\end{equation}
We are done. \endpf

\

\ni Now, let $\mu$ be a Borel $f$-invariant ergodic 
measure on $J(f)$ 
such that $h_{\mu }(f|\sg )>0$, where $h_{\mu }(f|\sg )$ denotes the 
relative entropy of $(f,\mu )$ with respect to $\sg $,\  
and let $T_\mu:L_\infty(\mu)\to L_\infty(\mu)$ be the Perron-Frobenius operator
associated to the measure $\mu$. It is defined by the formula
\begin{equation}\lab{4081605}
T_\mu g(x)=\sum_{y\in f^{-1}(x)}J_\mu^{-1}(y)g(y).
\end{equation}
Since $\L h=h$ everywhere throughout $J(f)$, using Lemma~6.9 from 
\cite{hiroki1}, we get that
\begin{equation}\lab{1081605}
\aligned
1
&=\int\1 d\mu
 =\int{\L h\over h}d\mu
 =\int T_\mu\lt({h\exp(\ov\psi)\over J_\mu^{-1}\cdot h\circ f}\rt)d\mu
 =\int{h\exp(\ov\psi)\over J_\mu^{-1}\cdot h\circ f} d\mu \\
&\ge 1+\int\log\lt({h\exp(\ov\psi)\over J_\mu^{-1}\cdot h\circ f}\rt)d\mu \\ 
& =1+\int\log h d\mu-\int\log(h\circ f)d\mu+\int\ov\psi d\mu
 +\int\log J_\mu d\mu \\
&=1+\int\psi d\mu-\P_b(\psi) + \hmu(f).
\endaligned
\end{equation}
Seeking contradiction suppose that $\P(\psi)>\P_b(\psi)$. By the 
Variational Principle there exists a Borel probability $f$-invariant 
ergodic measure $\mu$ such that
$$
\hmu(f)+\int\psi d\mu
>\P(\psi)-\min\{\P(\psi)-\P_b(\psi), \P(\psi)-\sup(\psi)-\log s\}.
$$
Note that, because of (\ref{6061805}), 
$\P(\psi)-\sup(\psi)-\log s>0$,
and therefore
$$
\aligned
\hmu(f)
& >\P(\psi)-(\P(\psi)-\sup(\psi)-\log s)-\int\psi d\mu \\ 
& =\log s+\sup(\psi)-\int\psi d\mu \\ 
& \ge  \log s,
\endaligned
$$which implies $h_{\mu }(f|\sg )>0.$
Hence, (\ref{1081605}) applies, and we can continue it to get that
$1\ge 1+\int\psi d\mu-\P_b(\psi)+\hmu(f)>1+\P(\psi)-(\P(\psi)-\P_b(\psi))
-\P_b(\psi)=1$. This contradiction along with (\ref{2081605}) and 
Lemma~\ref{l5061805} give the following.

\

\bprop\lab{p1081605}
$\P(\psi)=\P_b(\psi)$ and $\mu_\psi$ is an equilibrium state for
$\psi$. 
\eprop
\begin{lem}
\lab{rempbpressure}
For each point $x\in J(f)$, we have 
$\frac{1}{n}\log \L_{\psi }^{n}\1 (x)\rightarrow 
\P(\psi )$ as $n\rightarrow \infty .$ 
\end{lem}
{\sl Proof.}
By Proposition~\ref{p1081605}, $\L_{\psi }=(\exp (\P(\psi )))\cdot 
\L.$ Combining it with Lemma~\ref{l1061805} and Lemma~\ref{l1061805a}, 
we obtain that the statement of our lemma holds. \endpf

\

\ni Now suppose that $\mu$ is an arbitrary ergodic equilibrium state for 
$\psi$. Then $h_{\mu }(f)-\log s=\P(\psi )-\int \psi d\mu -\log s\geq 
\P(\psi )-\sup (\psi )-\log s>0.$ Hence 
$h_{\mu }(f|\sigma )>0.$ 
In view of the proposition above, the last component in (\ref{1081605}) 
is equal to $1$. Consequently, the only inequality in this 
formula becomes an equality, and so $\log\lt({h\exp(\ov\psi)
\over J_\mu^{-1}\cdot h\circ f}\rt)=0$ $\mu$-a.e. We thus get the 
following.

\

\blem\lab{l2081605}
If $\mu $ is an ergodic arbitrary equilibrium state for $\psi$, then $J_\mu
={h\circ f\over h}\exp(\P(\psi)-\psi )$ $\mu$-a.e.
\elem

\

\ni Now we shall prove the following.

\

\blem\lab{l2081605a}
If a Borel $f$-invariant probability measure $\mu $ satisfies 
$J_{\mu }=\frac{h\circ f}{h}\exp (\P(\psi )-\psi )\ \mu $-a.e., and 
if a Borel $f$-invariant probability measure $m$ satisfies 
$J_{m}=\frac{h\circ f}{h}\exp (\P(\psi )-\psi )$ everywhere,\ 
then $\mu $ is absolutely continuous with respect to $m.$ 
\elem
{\sl Proof.} 
It suffices to prove that for every $\e>0$ there exists $\d>0$ such that
if $g:J(f)\to (0,1]$ is a continuous function and $\int g dm\le\d$,
then $\int gd\mu\le \e$. Taking $q=q(\epsilon )\ge 1$ large enough, 
taking points $z_{j}$ ($j=1,\ldots\! ,u$) in 
$J(G)$ such that 
$\cup _{j=1}^{u}B(z_{j},R_{q}(z_{j}))\supset J(G)$, 
using 
(\ref{1061705}), Lemma~\ref{l1061805}, and redoing the considerations between 
(\ref{3081605}) and (\ref{3061805}) with the function $\1$ replaced by 
the function $g$, and 
$\| G_{j,\om }^{n} \1 \| _{\infty }$ replaced by 
$G_{j,\om }^{n}g(\tau , x)$, 
after applying 
(\ref{4081605}),
we get for every $1\le j\le u$, every $n=l\cdot q\ge 1\ (l\in \Bbb{N}$) and every $x\in\pi_2^{-1}
\(B(z_j,R_{q}(z_j))\)$, that
\begin{equation}\lab{7081605}
T_\mu^{qn}g(x)\le {\e\over 2}+CG_j^ng(x) \
\text{ and } \
T_{m}^{qn}g(x)\ge C^{-1}G_j^ng(x), 
\end{equation}
where $C=\overline{Q}_{\psi }/\underline{Q}_{\psi }$ is a positive constant which
does not depend on  
$l,j,x,g.$ 
We take a partition of the set $J(G)$ into mutually
disjoint sets $\{Y_j^i:1\le j\le u,\, 1\le i\le i(j)\}$ each of which
has a positive $m\circ \pi _{2}^{-1}$ measure and 
$\bu_{i=1}^{i(j)}Y_{j}^{i}\sbt
B(z_j,R_{q}(z_j))$ for all $j=1,2,\ld,u$. 
(Note that we can take such a partition, since 
the support of $m$ is equal to $J(f)$, which follows 
from the topological exactness of $f:J(f)\rightarrow J(f)$ and 
the condition 
``$J_{m }=\frac{h\circ f}{h}\exp (\P(\psi )-\psi )$ everywhere".)
Put 
$$\alpha =\sup\lt\{{\mu\circ \pi _{2}^{-1}\(Y_j^i\)
\over m\circ \pi _{2}^{-1}\(Y_j^i\)}:
  1\le j\le u,\, 1\le i\le i(j)\rt\}.$$
Since $g$ is everywhere positive, there exist $\zeta>0$ and $k\ge 1$
such that
\begin{equation}\lab{7081605a}
{1\over 2}\le {g(\th,y)\over g(\rho,z)}\le 2
\end{equation}
for all $(\th,y), (\rho,z)\in J(f)$ with $d(y,z)<\zeta$ and $\th|_k
=\rho|_k$. Fix now $n=l\cdot q\ge 1$ so large that 
for each $1\leq j\leq u$, 
$K_{G}(R'_{q}(z_{j}))^{4}\ka^{-1}\l^{n/2}C_{sa}^{1/2}
<\zeta$ (see item ($B_n$) of Corollary~\ref{c5042905}). Using then item
($B_n$) of this corollary and (\ref{7081605a}), the application of 
(\ref{6081605}) gives for all $n=l\cdot q\ge k$, all $\om\in \{1,\ld,s\}^n$, and 
all $\xi,x\in\pi_2^{-1}\(B(z_j,R_{q}(z_j))\)$, that
$$
\frac{1}{2}e^{-2H}G_{j,\om}^ng(x)\le G_{j,\om}^ng(\xi)\le 2e^{2H}G_{j,\om}^ng(x).
$$
Hence summing over all $\om\in \{1,\ld,s\}^n$, we get
\begin{equation}\lab{8081605}
\frac{1}{2}e^{-2H} G_j^ng(x)\le G_j^ng(\xi)\le 2e^{2H} G_j^ng(x)
\end{equation}
for all $n=l\cdot q\ge k$, all $1\le j\le u$ and all $\xi,x\in
\pi_2^{-1}\(B(z_j,R_{q}(z_j))\)$. 
Let $X_{j}^{i}:= \pi _{2}^{-1}(Y_{j}^{i}).$ 
Now applying (\ref{7081605}), 
(\ref{8081605}), and the definition of $\alpha $, we get
$$
\aligned
\int gd\mu
&=\int T_\mu^{qn}gd\mu
 \le {\e\over 2}+\sum_{j=1}^{u}C\sum_{i=1}^{i(j)}\int_{X_j^i}G_j^ngd\mu \\
&\le {\e\over 2}+\sum_{j=1}^{u}C\sum_{i=1}^{i(j)}
     \mu\(X_j^i\)\sup_{X_j^i}\(G_j^ng\) \\
&\le {\e\over 2}+\sum_{j=1}^{u}C\sum_{i=1}^{i(j)}
     2e^{2H}\mu\(X_j^i\)\inf_{X_j^i}\(G_j^ng\) \\
&\le {\e\over 2}+2e^{2H}\sum_{j=1}^{u}C\sum_{i=1}^{i(j)}{\mu\(X_j^i\)\over
     m\(X_j^i\)}\int_{X_j^i}G_j^ngdm \\
&\le {\e\over 2}+2e^{2H}\alpha \sum_{j=1}^{u}C^{2}\sum_{i=1}^{i(j)}
     \int_{X_j^i}T_{m}^{qn}gdm \\
&=   {\e\over 2}+2e^{2H}\alpha C^{2}\int T_{m}^{qn}g dm 
 =   {\e\over 2}+2e^{2H}\alpha C^{2}\int gdm.
\endaligned
$$
So, taking $\d=(4\alpha C^{2}e^{2H})^{-1}\e$, finishes the proof. \endpf

\

\blem\lab{l2081605b}
If $\mu $ is an ergodic equilibrium state, and if a Borel $f$-invariant 
probability measure $m$ satisfies 
$J_{m }=\frac{h\circ f}{h}\exp (\P(\psi )-\psi )$ everywhere,\ 
then $\mu $ is absolutely continuous with respect to $m.$ 
\elem
{\sl Proof.} 
By Lemma~\ref{l2081605} and Lemma~\ref{l2081605a}, we obtain the 
statement.
\endpf

\ 

\ni We conclude the paper with the following.

\

\bthm\lab{t4081605}
The measure $\mu_\psi$ forms a unique (ergodic) equilibrium state for $\psi$
and $m_\psi$ is a unique $\exp(P(\psi)-\psi)$-conformal measure.
\ethm
{\sl Proof.} For every $\exp(\P(\psi)-\psi)$-conformal measure $\nu$ let
$h=h^\nu$ be the density function obtained in Lemma~\ref{l3061805} with 
$m_\psi$ replaced by $\nu$. Note that by 
Remark~\ref{rnormal}, 
they differ by a positive multiplicative constant.
Suppose that there exist two $\exp(\P(\psi)-\psi)$-conformal measures
$m_1$ and $m_2$ that are not equivalent. Then there exists a Borel set $A\sbt 
J(f)$ such that $m_1(A)>0$ and $m_2(A)=0$. Hence, on the one hand, 
$h^{m_1}m_1(A)>0$
since, by Lemma~\ref{l3061805}, $h^{m_{1}}>0$ everywhere, and, on the other hand,
$h^{m_1}m_1(A)=0$ since each measure $\tau _{i}=h^{m_{i}}m_i\ (i=1,2)$ 
satisfies $J_{\tau _{i}}=\frac{h^{m_{i}}\circ f}{h^{m_{i}}}\exp (\P(\psi )-\psi )$ 
everywhere and 
hence these two measures are equivalent 
by Lemma~\ref{l2081605a}. So,
any two $\exp(\P(\psi)-\psi)$-conformal measures are equivalent. Now suppose
that there are two different ergodic equilibrium states $\mu_1$ and $\mu_2$
for $\psi$. Then they are mutually singular and there is a completely
invariant Borel set $A\sbt J(f)$ ($f(A)=A=f^{-1}(A)$) such that $\mu_1(A)=1$
and $\mu_2(A)=0$ (which implies that $\mu_2(J(f)\sms A)=1$). Then, by 
Lemma~\ref{l2081605b}, $m_\psi(A)>0$ and $m_\psi(J(f)\sms A)>0$. Hence,
both Borel probability measures $m_\psi^1$ and $m_\psi^2$ on $J(f)$, 
respectively defined by the formulas
$$
m_\psi^1(F)
={m_\psi(F\cap A)\over m_\psi(A)}\
\text{ and } \
m_\psi^2(F)
={m_\psi(F\cap(J(f)\sms A))\over m_\psi(J(f)\sms A)},
$$
are $\exp(\P(\psi)-\psi)$-conformal. Since they are singular, we get a
contradiction, and the uniqueness of equilibrium state for $\psi:J(f)\to\R$
is established. So, coming back to conformal measures, if $\nu_1$ and $\nu_2$ are
two $\exp(\P(\psi)-\psi)$-conformal measures, then $h^{\nu_1}\nu_1=
h^{\nu_2}\nu_2$. Since the ratio $h^{\nu_2}/h^{\nu_1}$ is constant, so is the
ratio $\nu_2/\nu_1$. Since, in addition, both these measures $\nu_2$ and 
$\nu_1$ are probabilistic, they are equal. We are done. \endpf

\

\noindent {\bf Proof of the main result Theorem~\ref{mainth}}: 
Combining Theorem~\ref{t4081605} and Lemma~\ref{rempbpressure}, 
the statement of the theorem follows.
\endpf

\end{document}